\documentclass[12pt]{article} 
\usepackage[mathlines]{lineno}
\usepackage[left=2cm,right=2cm,top=2cm,bottom=2cm]{geometry}
\usepackage{mathtools,amssymb,latexsym,amsthm,amsmath,amsfonts}   
\usepackage{xcolor}
\usepackage{hyperref}
\usepackage{authblk}
\usepackage{algorithm}
\usepackage{algpseudocode}
\usepackage{cite}
 \usepackage{url}
\usepackage{fancyhdr}
\usepackage{amsthm}
\usepackage[all]{xy}
\usepackage{fontenc}
\usepackage{listings}

\newtheorem{theorem}{Theorem}[section]
  \newtheorem{proposition}{Proposition}[section]
               
               \newtheorem{example}{Example}[section]

               \newtheorem{remark}{Remark}[section]
               
               \def\pf{\par\noindent {\em Proof.}~\par\noindent}
               \def\qed{~\hfill{$\square$}\pagebreak[1]\par\medskip\par}

\newcommand{\B}{{\cal B}}
\newcommand{\R}{{\mathbb R}}

\newcommand{\I}{{\cal I}}

\newcommand{\ux}{\underline{x}}
\newcommand{\pD}{{^\psi\!\pux}}
\newcommand{\hD}{{^\varphi\!\pux}}

\newcommand{\wH}{{^{\varphi,\psi}\!\Psi}}

\newcommand{\pux}{\partial_{\ux}}

\newcommand{\cP}{{\cal P}}

\begin{document}

\title{Revisiting Fischer decompositions by inframonogenic functions}
\author[1]{Daniel Alfonso Santiesteban}
\author[1]{Ricardo Abreu Blaya}	
\author[2]{Juan Bory Reyes}
\author[3]{Baruch Schneider\footnote{Corresponding author}}

\affil[1]{Facultad de Matem\'aticas, Universidad Aut\'onoma de Guerrero, Chilpancingo de los Bravo, Mexico\\\href{mailto:danielalfonso950105@gmail.com}{danielalfonso950105@gmail.com}, \href{mailto:rabreublaya@yahoo.es}{rabreublaya@yahoo.es}}
\affil[2]{ESIME, Instituto Polit\'ecnico Nacional, CDMX, M\'exico\\\href{mailto:juanboryreyes@yahoo.com}{juanboryreyes@yahoo.com}}
\affil[3]{Department of Mathematics, University of Ostrava, Czech Republic\\\href{mailto:baruch.schneider@osu.cz}{baruch.schneider@osu.cz}}

\maketitle
\begin{abstract}
We explore Fischer decompositions from the perspective of inframonogenic functions in Clifford analysis, a non-commutative version of harmonic functions in Euclidean spaces. We study decomposition formulas for polynomial spaces involving inframonogenic and harmonic components, providing fresh insights into their algebraic structures and further important applications. 
\end{abstract}

\vspace{0.2cm}

\small{
}
\section{Introduction}
Clifford analysis is a higher-dimensional generalization of classical complex function theory in which Clifford algebras provide the underlying algebraic framework. It combines tools from harmonic analysis, differential equations, and algebra to study Clifford algebra-valued functions defined on Euclidean spaces. A central role in this function theory is played by the Dirac operator, which serves as a higher-dimensional analogue of the Cauchy--Riemann operator and gives rise to the notion of monogenic functions through the associated generalized Cauchy--Riemann equations \cite{BDS}. This framework provides a unified setting for investigating harmonic, monogenic, and more general classes of solutions of elliptic systems, as well as their algebraic and analytic properties.

Let $\mathcal{P}(\mathbb{R}^{m},\R_{0,m})$ denote the space of $\R_{0,m}$-valued polynomials on $\mathbb{R}^{m}$, where $\R_{0,m}$ is the $2^m$-dimensional real Clifford algebra constructed over the orthonormal basis $\{e_1,...,e_m\}$ of the Euclidean space $\R^m$ and determined by relations $e_je_i+e_ie_j=-2\delta_{ji}$. Let
\[
\ux=\sum_{j=1}^{m}x_j e_j
\]
denote the vector variable in the Euclidean space $\mathbb{R}^m$. The Euclidean Dirac operator is defined by
\[
\pux=\sum_{j=1}^{m}e_j\partial_{x_j}.
\]
In particular, the equation
\[
\pux f=0
\]
defines the class of left-monogenic functions, which may be regarded as a higher-dimensional analogue of holomorphic functions. Every monogenic function is harmonic, as a consequence of the factorization of the $m$-dimensional Laplacian by the Dirac operator:
\[\Delta_m=-\pux^2.\]
The corresponding spaces of homogeneous monogenic polynomials provide the basic building blocks for several important decompositions in Clifford analysis. For more information we refer the reader to \cite{BDGS,BDShS,DSS}.

Fischer decompositions constitute one of the fundamental algebraic tools in Clifford analysis and provide a natural framework for understanding the structure of polynomial solutions of differential equations associated with Dirac-type operators. Their origins can be traced back to the classical Fischer decomposition in harmonic analysis, where homogeneous polynomials are decomposed into orthogonal components involving powers of the Euclidean norm and harmonic polynomials. Let $q(\ux)$ be a homogeneous $\R_{0,m}$-valued polynomial in $\mathbb{R}^m$. Every homogeneous polynomial $P_k(\ux)$ of degree $k$ admits a unique decomposition of the form
\[
P_k(\ux)=Q_k(\ux)+q(\ux)R_{k-d}(\ux),
\]
where $d=\deg q$, $Q_k$ is a homogeneous polynomial of degree $k$ satisfying
\[
q(\partial_{\ux})Q_k=0,
\]
and $R_{k-d}$ is a homogeneous polynomial of degree $k-d$ (with $R_{k-d}=0$ if $k<d$). Here, $q(\partial_{\ux})$ denotes the constant-coefficient differential operator associated with $q$, obtained by replacing each variable $x_j$ in $q(\ux)$ with the corresponding partial derivative $\partial_{x_j}$. This decomposition is commonly referred to as the \emph{Fischer decomposition}, following the work of Ernst Fischer in 1917 (see \cite{EF}).

A fundamental example is obtained by taking
\[
q(\ux)=|\ux|^2=\sum_{j=1}^mx_j^2.
\]
In this case,
$q(\partial_{\ux})=\Delta_m$,
 and the condition $q(\partial_{\ux})Q_k=0$ reduces to the harmonicity condition $\Delta_m Q_k=0$. Consequently, the classical Fischer decomposition yields a decomposition of polynomials into components generated by powers of $|\ux|^2$ and homogeneous harmonic polynomials, namely,
 \begin{equation}\label{harmonicFischer}
\mathcal{P}(\mathbb{R}^m,\R_{0,m})
=
\bigoplus_{k=0}^\infty\bigoplus_{p=0}^{\infty}
|\ux|^{2p}\mathcal{H}(k), 
 \end{equation}
where $\mathcal{H}(k)$ denotes the subspace of $\mathcal{P}(\mathbb{R}^m,\R_{0,m})$ consisting of harmonic homogeneous polynomials of degree $k$. This decomposition, and its various extensions, plays a fundamental role in harmonic ana\-lysis, representation theory, partial differential equations, and mathematical physics. In particular, its formulation in terms of solid spherical harmonics provides a powerful algebraic and analytic tool for the study of polynomial solutions of elliptic equations and for the construction of orthogonal expansions \cite{BDShS,brackx2,DSS, fischer4}. Furthermore, the Fischer decomposition can be interpreted as an irreducible decomposition of $\mathcal{P}(\R^m,\R_{0,m})$ under the action of $\mathfrak{osp}(m|2n)$. More precisely, $\mathfrak{osp}(m|2n)$ and $\mathfrak{sl}(2)$ constitute a Howe dual pair acting on the polynomial superalgebra over $\mathbb{R}^{m|2n}$, except in the singular cases where the superdimension $M=m-2n$ belongs to $-2\mathbb{N}_0$; see \cite{coul}. From now on, the terminology ``$k$-homogeneous polynomial'' is used to refer to homogeneous polynomials of degree $k$.

Clifford analysis allows one to obtain a refinement of the decomposition \eqref{harmonicFischer} by means of monogenic functions. If one takes $q(\ux)=\ux$, then $q(\pux)=\pux$, and $Q_k$ becomes a monogenic $k$-homogeneous polynomial. This leads to the Fischer decomposition
\begin{equation}\label{farm2}
\mathcal{P}(\mathbb{R}^m,\mathbb{R}_{0,m})
=
\bigoplus_{k=0}^\infty
\bigoplus_{p=0}^\infty
\ux^p \mathcal{M}(k),
\end{equation}
where $\mathcal{M}(k)\subset\mathcal{P}(\R^m,\R_{0,m})$ is the subspace of left-monogenic $k$-homogeneous polynomials. The Fischer decomposition \eqref{farm2}, together with the Cauchy--Kovalevskaya extension, constitutes one of the fundamental tools for constructing orthonormal bases of spaces of monogenic polynomials
\cite{brackx1,KV}.

Of particular relevance to the present work is the paper \cite{MPS2}, where the authors established a new Fischer decomposition for arbitrary homogeneous polynomials in terms of solutions to the sandwich equation
\[
\pux f\pux = 0,
\]
which are referred to as inframonogenic functions; see, for instance, \cite{ABMM,MPS1,MAB1,MAB2,MAB3,MMA}. More precisely, they obtained the complete Fischer decomposition
\begin{equation}\label{fischerinfra}
\cP(\R^m,\R_{0,m})=\bigoplus_{k=0}^{\infty}\bigoplus_{p=0}^{\infty}\ux^p\I(k)\ux^p,
\end{equation}
where $\I(k)$ denotes the space of inframonogenic $k$-homogeneous polynomials in $\R^m$. More recently, an improvement of this result was obtained in \cite{La} by means of the so-called $H$-action of the Pin group on Clifford algebra-valued polynomials.

On the other hand, as demonstrated in several works, significant developments in Clifford analysis have been achieved by replacing the standard basis $\{e_1,e_2,\dots,e_m\}$ with an arbitrary orthonormal basis of $\R^m$
\[
\{\varphi_1,\varphi_2,\dots,\varphi_m\},
\]
known as a structural set \cite{SV}, and by considering the corresponding Dirac operator
\[
\hD=\sum_{j=1}^m\varphi_j\partial_{x_j}.
\]
In the setting of two different structural sets, a mixed Fischer decomposition was established in \cite{BDGS} in terms of $(\varphi,\psi)$-monogenic functions, namely, functions satisfying the overdetermined system
\[
\hD f=0,\qquad \pD f=0.
\]
The consideration of arbitrary structural sets provides a natural framework for studying the more general second-order partial differential equation
\[
\hD f\pD=0,
\]
whose solutions are referred to as $(\varphi,\psi)$-inframonogenic functions \cite{A1,APA1}. This greater flexibility opens up new perspectives in several directions of research, including the mapping properties of $\Pi$-operators, alternative Kolosov--Muskhelishvili formulas, additive representations of contragenic polynomials, and geometric conformal mappings (see, for example, \cite{KM1,ABGU,GN1,GN2}). 

In \cite{FischerDAS}, the preceding framework was further developed, leading to the following complete Fischer decomposition in terms of $(\varphi,\psi)$-inframonogenic polynomials:
\begin{equation}\label{fischerinfrageneral}
\cP(\R^m,\R_{0,m})=\bigoplus_{k=0}^{\infty}\bigoplus_{p=0}^{\infty}\ux_\varphi^p\I_{\varphi,\psi}(k)\ux_\psi^p,
\end{equation}
where $\ux_\varphi=\sum_{j=1}^m\varphi_jx_j$, $\ux_\psi=\sum_{j=1}^m\psi_jx_j$ and $\I_{\varphi,\psi}(k)$ denotes the space of $(\varphi,\psi)$-inframonogenic $k$-homogeneous polynomials in $\R^m$.

In the present paper, we investigate several subtleties and interesting consequences arising from the Fischer decompositions established in \cite{FischerDAS} and \cite{MPS2}, particularly from the decomposition in terms of inframonogenic polynomials  \eqref{fischerinfra}  and its more general counterpart involving $(\varphi,\psi)$-inframonogenic polynomials  \eqref{fischerinfrageneral}. We show that these decompositions exhibit essential differences from the Fischer decomposition in terms of biharmonic polynomials. 

\section{Preliminaries}

We begin by recalling some basic notions from Clifford analysis that will be used throughout the paper. Let $\R_{0,m}$ denote the real Clifford algebra associated with the Euclidean space $\R^m$, generated by the orthonormal basis $\{e_1,\ldots,e_m\}$ and subject to the relations
\[
e_j e_i+e_i e_j=-2\delta_{ji},
\qquad j,i=1,\ldots,m.
\]
Thus, $\R_{0,m}$ is a $2^m$-dimensional real associative algebra. If
$A=\{j_1,\ldots,j_n\}\subseteq\{1,\ldots,m\}$, with
$j_1<\cdots<j_n$, we write
\[
e_A=e_{j_1}\cdots e_{j_n},
\]
and set $e_{\emptyset}=1$, where $1$ denotes the identity element of
$\R_{0,m}$. Every element $a\in\R_{0,m}$ admits a unique representation
\[
a=\sum_A a_A e_A,
\qquad a_A\in\R.
\]

The Clifford algebra is naturally graded. More precisely,
\[
\R_{0,m}=\bigoplus_{k=0}^m \R_{0,m}^{(k)},
\]
where $\R_{0,m}^{(k)}$ denotes the subspace of $k$-vectors,
\[
\R_{0,m}^{(k)}
=
\left\{
\sum_{|A|=k}a_Ae_A:\ a_A\in\R
\right\}.
\]
Accordingly, every $a\in\R_{0,m}$ can be uniquely decomposed into its $k$-vector components as
\[
a=\sum_{k=0}^m [a]_k,
\qquad [a]_k\in\R_{0,m}^{(k)},
\]
where $[a]_k$ denotes the projection of $a$ onto the subspace of
$k$-vectors. In particular, $\R_{0,m}^{(0)}$, $\R_{0,m}^{(1)}$, $\R_{0,m}^{(0)}\oplus\R_{0,m}^{(1)}$, and $\R_{0,m}^{(m)}$ are referred to, respectively, as the scalar, vector, paravector, and pseudoscalar subspaces of $\R_{0,m}$. The associative Clifford algebra $\R_{0,m}$ admits the decomposition into the direct sum of its even and odd subspaces,
\[
\R_{0,m}=\R_{0,m}^{+}\oplus\R_{0,m}^{-},
\]
where $\R_{0,m}^{+}$ and $\R_{0,m}^{-}$ are $2^{m-1}$-dimensional and consist, respectively, of the even and odd multivectors. Moreover, $\R_{0,m}^{+}$ is a subalgebra of $\R_{0,m}$, commonly referred to as the even subalgebra. Consequently, every Clifford number $a\in\R_{0,m}$ admits a unique decomposition
\[
a=a_{+}+a_{-},\qquad a_{+}\in\R_{0,m}^{+},\quad a_{-}\in\R_{0,m}^{-},
\]
where $a_{+}$ and $a_{-}$ are called the even and odd parts of $a$, respectively.

We shall also use the standard Clifford conjugation and reversion, which are the
anti-involution determined by
\[
\overline{e_j}=-e_j,\quad \widehat{e_j}=e_j,
\qquad j=1,\ldots,m,
\]
and
\[
\overline{ab}=\overline{b}\,\overline{a},\quad \widehat{ab}=\widehat{b}\,\widehat{a},
\qquad a,b\in\R_{0,m}.
\]
The Euclidean norm on $\R_{0,m}$ is defined by
\[
|a|^2=[a\overline{a}]_0.
\]
Under the canonical identification of $\R^m$ with the vector subspace
$\R_{0,m}^{(1)}$ of the Clifford algebra, this norm coincides with the
usual Euclidean norm; namely,
\[
|\ux|^2=x_1^2+\cdots+x_m^2,
\qquad \ux=x_1e_1+\cdots+x_me_m\in\R^m.
\]

Throughout the paper, functions will be defined on open subsets of
$\R^m$ and take values in $\R_{0,m}$. Any such function can be written
componentwise as
\[
f(x)=\sum_A f_A(x)e_A,
\]
where the coefficient functions $f_A$ are real-valued. Consequently, properties such as continuity, differentiability, and
integrability are understood componentwise.

The Euclidean Dirac operator is the first-order differential operator
\[
\pux
=
\sum_{j=1}^m e_j\partial_{x_j}
=
e_1\partial_{x_1}+\cdots+
e_m\partial_{x_m}.
\]
It satisfies the fundamental factorization
\[
\pux^2=-\Delta_m,
\]
where
\[
\Delta_m=\sum_{j=1}^m\partial_{x_j}^2
\]
is the Laplace operator in $\R^m$.

Let $\Omega\subset\R^m$ be open. An $\R_{0,m}$-valued function
$f\in C^1(\Omega)$ is called \emph{left-monogenic} in $\Omega$ if
\[
\pux f=0,
\]
and \emph{right-monogenic} in $\Omega$ if
\[
f\pux=0.
\]
Monogenic functions constitute one of the fundamental classes of
solutions in Clifford analysis; see, for example,
\cite{BDS,GHS,DSS}.

More generally, let
\[
\varphi=\{\varphi_1,\ldots,\varphi_m\}
\]
be an orthonormal basis of $\R^m$, regarded as a structural set in
$\R_{0,m}$. The corresponding Dirac operator is defined by
\[
\hD
=
\sum_{j=1}^m
\varphi_j\partial_{x_j}
=
\varphi_1\partial_{x_1}
+\cdots+
\varphi_m\partial_{x_m}.
\]
Likewise, for another structural set
\[
\psi=\{\psi_1,\ldots,\psi_m\},
\]
we consider the associated Dirac operator
\[
\pD
=
\sum_{j=1}^m
\psi_j\partial_{x_j}.
\]
Since both $\varphi$ and $\psi$ are orthonormal bases, the corresponding
Dirac operators satisfy
\[
\hD\!\!^2=\pD\!\!^2=-\Delta_m.
\]
The use of arbitrary structural sets provides a natural extension of
the classical setting based on the standard basis and has proved useful
in several areas of Clifford analysis
\cite{GN1,GN2,SV}. For
$\ux=\sum_{i=1}^m x_i e_i\in\R^m$, we introduce the Clifford-vector
notation
\[
\ux_\varphi=\sum_{i=1}^m x_i\varphi_i,
\qquad
\ux_\psi=\sum_{i=1}^m x_i\psi_i.
\]

A function $f\in C^1(\Omega,\R_{0,m})$ satisfying
\[
\hD f=0
\]
is called left-$\varphi$-monogenic, whereas a function satisfying
\[
f\hD=0
\]
is called right-$\varphi$-monogenic. Analogous terminology is
used for the structural set $\psi$ and the operator $\pD$.

In the present paper, however, our main interest lies in the
second-order equation obtained by coupling two, in general distinct,
structural sets. To this end, for an open set $\Omega\subset\R^m$, we
define
\[
\I_{\varphi,\psi}(\Omega)
=
\left\{
f\in C^2(\Omega,\R_{0,m}):
\hD f\pD=0
\right\}.
\]
The elements of $\I_{\varphi,\psi}(\Omega)$ will be referred to as
\emph{$(\varphi,\psi)$-inframonogenic functions}.

This definition contains the classical notion of inframonogenicity as
a particular case. Indeed, when
\[
\varphi=\psi
=
\{e_1,\ldots,e_m\},
\]
we have $\hD=\pD=\pux$, and consequently
$\I_{\varphi,\psi}(\Omega)
=\I(\Omega)$, the class of standard inframonogenic functions. 
Thus, $(\varphi,\psi)$-inframonogenic functions provide a natural
two-structural-set generalization of inframonogenic functions, which
were introduced and studied in \cite{MPS2,MPS1} (see also
\cite{A1}).

We shall denote the space of $\R_{0,m}$-valued $k$-homogeneous
polynomials by $\cP(k)$. Throughout the remainder of the paper, we shall use the following results interchangeably.
\begin{proposition}\cite[Theorem 4.3, p. 8]{FischerDAS}\label{Prop1}
Let $P_k\in\mathcal{P}(k)$ and let $s$ be a positive integer, then we have
\begin{equation}
\hD (\ux_{\varphi}^s P_k)=g_{s,k}\ux^{s-1}P_k+(-1)^s\ux_{\varphi}^s(\hD P_k),
\end{equation}
where $g_{2n,k}=-2n$ and $g_{2n+1,k}=-(2n+2k+m)$. Similarly,
\begin{equation}
( P_k \ux_{\varphi}^s)\hD=g_{s,k}P_k \ux_{\varphi}^{s-1}+(-1)^s( P_k\hD)\ux_\varphi^s.
\end{equation}
\end{proposition}
Consider the linear operator
\begin{equation}
\wH_1 (f)=\sum_{j=1}^m\varphi_j f\psi_j.
\end{equation}
In the following relations involving $\pD$, $\hD$, and $\wH_1$, which can be found in papers \cite{APA1} and \cite{DASGLNS}, we tacitly assume that all required differentiations are well defined.
\begin{proposition}\label{Prop2}
Let be $f:\R^m\to\R_{0,m}$, then
\begin{itemize}
\item[(i)] $\hD[\wH_1(f)]=-2f\pD-\wH_1(\hD f)$,\quad   $[\wH_1(f)]\pD=-2\hD f-\wH_1(f\pD)$,
\item[(ii)] $\hD[\wH_1(f)]\pD=\wH_1(\hD f\pD)$,\quad $\Delta_m\wH_1(f)=\wH_1(\Delta_m f)$,
\item[(iii)] $\hD[f\ux_\psi]=(\hD f)\ux_\psi+\wH_1(f)$,\quad $[\ux_\varphi f]\pD=\ux_\varphi (f\pD)+\wH_1(f)$,
\item[(iv)]$\hD[f\ux^2]=(\hD f)\ux^2-2\ux_\varphi f$,\quad $[\ux^2 f]\pD=\ux^2(f\pD)-2f\ux_\psi$.
\end{itemize}
\end{proposition}
\section{Operators $\wH_n$}
Let $\varphi=\{\varphi_1,...,\varphi_m\}$ and $\psi=\{\psi_1,...,\psi_m\}$ be two arbitrary structural sets in $\R^m$. In \cite{APA1}, the following operators $\wH_n:\R_{0,m}\to\R_{0,m}$ were defined for $n=1,...,m$:
$$\wH_n(a)=\sum_{|A|=n}\varphi_A a\widehat{\psi_A},$$
where $\varphi_A=\varphi_{j_1}\cdots\varphi_{j_n}$ and $\psi_A=\psi_{j_1}\cdots\psi_{j_n}$. These operators play an important role in establishing necessary and sufficient conditions for inframonogenicity. In the same paper \cite{APA1}, the authors proved that the operator $\wH_1$ is a bijective $\R$-linear mapping on $\R_{0,m}$ when $m$ is odd. However, this property does not hold in even dimensions, where injectivity fails. Indeed, if we consider in $\R_{0,2}$ the structural sets $\varphi=\{e_1,e_2\}$ and $\psi=\{e_2,e_1\}$, then a direct computation shows that $\wH_1(1)=0$ and $\wH_1(e_1e_2)=0$. Nevertheless, in \cite{DASABory}, the authors proved the following proposition, which provides a certain bijectivity property of this operator on specific $k$-vector subspaces, making use of the transition matrix $\mathcal{M}_{\varphi,\psi}$ from the basis $\varphi$ to $\psi$.
\begin{proposition}{\cite[Proposition 3.1, p. 8]{DASABory}}\label{Prop1das}
Let $\varphi$ and $\psi$ be two arbitrary structural sets and $m$ be an even number.  If
$m\equiv 0{\pmod{4}}$ and $\det(\mathcal{M}_{\varphi,\psi})=1$, either $m\equiv 2{\pmod{4}}$ and $\det(\mathcal{M}_{\varphi,\psi})=-1$, then $\wH_1$ is a bijective $\R$-linear mapping on $\R_{0,m}^-$. Whereas, if  $m\equiv 2{\pmod{4}}$ and $\det(\mathcal{M}_{\varphi,\psi})=1$, either $m\equiv 0{\pmod{4}}$ and $\det(\mathcal{M}_{\varphi,\psi})=-1$, then $\wH_1$ is a bijective $\R$-linear mapping on $\R_{0,m}^+$.
\end{proposition}
A direct computation yields
\begin{align}\label{rel2}
\varphi_1\cdots\varphi_m
=
\det(\mathcal{M}_{\varphi,\psi})\psi_1\cdots\psi_m.
\end{align}
On the other hand, let $a_k\in\R_{0,m}^{(k)}$ and $s\in\R_{0,m}^{(m)}$. It is well known that, when $m$ is even,
\begin{align}\label{rel3}
sa_k&=a_ks,\qquad &&\text{if } k \text{ is even},\nonumber\\
sa_k&=-a_ks,\qquad &&\text{if } k \text{ is odd}.
\end{align}
In odd dimensions, the pseudoscalar commutes with every Clifford number.
The following recursive formula, proved in \cite{APA1}, will be useful in what follows:
\begin{equation}\label{recformula}
(m-n+1)\wH_{n-1}(f)+(n+1)\wH_{n+1}(f)=\wH_1(\wH_n(f)),\quad n=1,...,m-1.
\end{equation}

Consider the equation 
\begin{equation}
\wH_1(f)=g,
\end{equation}
where $g\in\R_{0,m}$. We are interested in determining $f$ when $g$ is known, whenever possible. To this end, we will conveniently use the facts established above, together with formula \eqref{recformula}. It is readily seen that
\begin{align*}
mf+2\wH_2(f)=\wH_1(g),
\end{align*}
and, thus,
\begin{align*}
m\wH_1(f)+2\wH_1(\wH_2(f))=\wH_1(\wH_1(g)).
\end{align*}
Hence,
$$mg+2(m-1)g+6\wH_3(f)=\wH_1^2(g),$$
and, therefore,
$$\wH_3(f)=\frac{1}{6}[\wH_1^2(g)-(3m-2)g].$$
By repeatedly applying formula \eqref{recformula}, we obtain
\begin{align*}
&(m-2)\wH_2(f)+4\wH_4(f)=\frac{1}{6}[\wH_1^3(g)-(3m-2)\wH_1(g)],\\
&\frac{(m-2)}{2}[\wH_1(g)-mf]+4\wH_4(f)=\frac{1}{6}[\wH_1^3(g)-(3m-2)\wH_1(g)],\\
&-\frac{m(m-2)}{2}f+4\wH_4(f)=\frac{1}{6}[\wH_1^3(g)-(6m-8)\wH_1(g)],\\
&-\frac{m(m-2)}{2}g+4(m-3)\wH_3(f)+20\wH_5(f)=\frac{1}{6}[\wH_1^4(g)-(6m-8)\wH_1^2(g)],\\
&20\wH_5(f)=\frac{1}{6}[\wH_1^4(g)-(6m-8)\wH_1^2(g)]+\frac{m(m-2)}{2}g-\frac{2(m-3)}{3}[\wH_1^2(g)-(3m-2)g],\\
&\wH_5(f)=\frac{1}{5!}[\wH_1^4(g)-(10m-20)\wH_1^2(g)+(15m^2-50m+24)g].
\end{align*}
In summary, we can continue and get the following
\begin{align*}
\wH_2(f)&=\frac{1}{2!}[\wH_1(g)-mf],\\
\wH_3(f)&=\frac{1}{3!}[\wH_1^2(g)-(3m-2)g],\\
\wH_4(f)&=\frac{1}{4!}[\wH_1^3(g)-2(3m-4)\wH_1(g)+3m(m-2)f],\\
\wH_5(f)&=\frac{1}{5!}[\wH_1^4(g)-10(m-2)\wH_1^2(g)+(15m^2-50m+24)g],\\
\wH_6(f)&=\frac{1}{6!}[\wH_1^5(g)-5(3m-8)\wH_1^3(g)+(45m^2-210m+184)\wH_1(g)-15m(m-2)(m-4)f],\\
\wH_7(f)&=\frac{1}{7!}[\wH_1^6(g)-7(3m-10)\wH_1^4(g)+(105m^2-630m+784)\wH_1^2(g)\\&\quad\quad\quad-(105m^3-840m^2+1764m-720)g],\\
\wH_8(f)&=\frac{1}{8!}[\wH_1^7(g)+(112-28m)\wH_1^5	+(210m^2-1260m+2464)\wH_1^3(g)\\
&\quad\quad+(-420m^3+4200m^2-11872m+8448)\wH_1(g)+105m(m-2)(m-4)(m-6)f],
\end{align*}
\begin{align*}
\wH_9(f)&=\frac{1}{9!}[\wH_1^8(g)+(168-36m)\wH_1^6(g)+(378m^2-1260m+6384)\wH_1^4(g)\\
&\quad\quad +(-1260m^3+15120m^2-53424m+52352)\wH_1^2(g)\\&\quad\quad+(945m^4-13860m^3+65772m^2-109584m+40320)g],\\
\wH_{10}(f)&=\frac{1}{10!}[\wH_1^9(g)+(240-45m)\wH_1^7(g)+(630m^2-3150m+14448)\wH_1^5(g)\\
&\quad\quad\quad+(-3150m^3+44100m^2-186480m+229760)\wH_1^3(g)\\
&\quad\quad\quad+(4725m^4-81900m^3+475020m^2-1040400m+648576)\wH_1(g)\\
&\quad\quad\quad-945m(m-2)(m-4)(m-6)(m-8)f].
\end{align*}

Now we will proceed to find a general formula. 

Let $D:=\wH_1$. Then,
\begin{equation}
\wH_n(f)=\frac{1}{n!}P_n(D)f,
\end{equation}
where 
\begin{equation}
P_0(x)=1,\quad P_1(x)=x,
\end{equation}
and
\begin{equation}
P_{n+1}(x)=xP_n(x)-n(m-n+1)P_{n-1}(x).
\end{equation}
Note that the first polynomials are
\begin{align*}
P_2(x)&=x^2-m,\\
P_3(x)&=x^3-(3m-2)x,\\
P_4(x)&=x^4-2(3m-4)x^2+3m(m-2),\\
P_5(x)&=x^5-10(m-2)x^3+(15m^2-50m+24)x.
\end{align*}
Let's write
\begin{equation}
P_n(x)=\sum_{j=0}^{\lfloor n/2\rfloor}(-1)^jc_{n,j}(m)x^{n-2j}.
\end{equation}
Therefore, the coefficients are explicitly given by
\begin{equation}
c_{n,j}(m)=\sum_{\overset{1\leq r_1<...<r_j\leq n-1}{r_{l+1}-r_l\geq 2}}\prod_{l=1}^jr_l(m-r_l+1),
\end{equation}
under the convention
$$c_{n,0}(m)=1.$$
Thus,
\begin{equation}
P_n(x)=\sum_{j=0}^{\lfloor n/2\rfloor}(-1)^jx^{n-2j}\sum_{\overset{1\leq r_1<...<r_j\leq n-1}{r_{l+1}-r_l\geq 2}}\prod_{l=1}^jr_l(m-r_l+1).
\end{equation}
The recurrence relation can be interpreted as the expansion of a tridiagonal determinant:
\begin{align*}
P_n(x)=\det\begin{pmatrix}
x&-1&0&\cdots&0\\
-m&x&-1&\cdots&0\\
0&-2(m-1)&x&\ddots&\vdots\\
\vdots&\ddots&\ddots&\ddots&-1\\
0&\cdots&0&-(n-1)(m-n+2)&x 
\end{pmatrix}.
\end{align*}
Upon expanding this determinant, each term in which $x$ is omitted at $j$ positions corresponds to selecting $j$ pairs of non-adjacent positions. This gives rise precisely to the sum
$$r_1<r_2<...<r_j,\quad r_{l+1}\geq r_l+2.$$
For example, for $n=6$,
$$c_{6,1}(m)=\sum_{r=1}^5r(m-r+1)=15m-40$$
and
$$c_{6,2}(m)=\sum_{\overset{1\leq r<s\leq 5}{s\geq r+2}}r(m-r+1)s(m-s+1)=45m^2-210m+184,$$
whereas
$$c_{6,3}(m)=15m(m-2)(m-4).$$
So,
$$P_6(x)=x^6-(15m-40)x^4+(45m^2-210m+184)x^2-15m(m-2)(m-4),$$
obtaining the expression we already know
\begin{align*}
\wH_6(f)&=\frac{1}{6!}[\wH_1^6(f)-5(3m-8)\wH_1^4(f)+(45m^2-210m+184)\wH_1^2(f)-15m(m-2)(m-4)f]\\
&=\frac{1}{6!}[\wH_1^5(g)-5(3m-8)\wH_1^3(g)+(45m^2-210m+184)\wH_1(g)-15m(m-2)(m-4)f].
\end{align*}
The algorithm can be implemented in MATLAB, Python, Julia, or other programming languages. We propose the following version in LaTeX pseudocode, which can then be translated into MATLAB with virtually no modifications.
\begin{algorithm}
\caption{Computation of $\wH_n(f)$}
\label{alg:Fn}
\begin{algorithmic}[1]

\Require $m,n\in\mathbb{N}$ with $n\geq 1$, the function $f$,
and the operator $\mathcal{F}_1$

\Ensure $\wH_n(f)$

\State $\wH_0\gets f$
\State $\wH_1\gets \mathcal{F}_1(\wH_0)$

\If{$n=0$}
    \State \Return $\wH_0$
\EndIf

\If{$n=1$}
    \State \Return $\wH_1$
\EndIf

\For{$k=1,\ldots,n-1$}
    \State
    $\wH_{k+1}\gets
    \dfrac{\mathcal{F}_1(\wH_k)-(m-k+1)\wH_{k-1}}{k+1}$
\EndFor

\State \Return $\wH_n$

\end{algorithmic}
\end{algorithm}

If $m$ is odd
\begin{equation}
\wH_m(f)=\frac{1}{m!}\sum_{j=0}^{(m-1)/2}(-1)^jc_{m,j}(m)\wH_1^{m-2j-1}(g),
\end{equation}
whence
\begin{equation}\label{form20}
f=\frac{\wH_m(1)}{m!}\sum_{j=0}^{(m-1)/2}(-1)^jc_{m,j}(m)\wH_1^{m-2j-1}(g).
\end{equation}

If $m$ is even
\begin{equation}\label{for21}
\wH_m(f)=(-1)^{m/2}f+\frac{1}{m!}\sum_{j=0}^{m/2-1}(-1)^jc_{m,j}(m)\wH_1^{m-2j-1}(g),
\end{equation}
whence
\begin{equation}\label{for22}
(1-(-1)^{m/2}\wH_m(1))f_+-(1+(-1)^{m/2}\wH_m(1))f_-=\frac{\wH_m(1)}{m!}\sum_{j=0}^{m/2-1}(-1)^jc_{m,j}(m)\wH_1^{m-2j-1}(g).
\end{equation}
\begin{remark}
Note that, when $m$ is odd, formula \eqref{form20} ensures the bijectivity of the operator $\wH_1$. In fact, if $g=0$, then $f=0$.   When m is even, the situation is different, as evidenced by Proposition \ref{Prop1das} and now by formula \eqref{for22}. Let us examine this case. First, observe that $\wH_m(1)=(-1)^m\det(\mathcal{M}_{\varphi,\psi})$. If
$m\equiv 0{\pmod{4}}$ and $\det(\mathcal{M}_{\varphi,\psi})=1$, either $m\equiv 2{\pmod{4}}$ and $\det(\mathcal{M}_{\varphi,\psi})=-1$, then
\begin{equation*}
f_-=-\frac{\wH_m(1)}{2\cdot m!}\sum_{j=0}^{m/2-1}(-1)^jc_{m,j}(m)\wH_1^{m-2j-1}(g),
\end{equation*}
hence, in general, the operator $\wH_1$ can only be bijective on $\R_{0,m}^-$. If  $m\equiv 2{\pmod{4}}$ and $\det(\mathcal{M}_{\varphi,\psi})=1$, either $m\equiv 0{\pmod{4}}$ and $\det(\mathcal{M}_{\varphi,\psi})=-1$, then
\begin{equation*}
f_+=\frac{\wH_m(1)}{2\cdot m!}\sum_{j=0}^{m/2-1}(-1)^jc_{m,j}(m)\wH_1^{m-2j-1}(g),
\end{equation*}
and the operator $\wH_1$ can only be bijective on $\R_{0,m}^+$. Note also that, when $g=0$, relation \eqref{for21} reduces to relation (3.1) presented in \cite{DASABory}.
\end{remark}
\begin{remark}\label{Remark3.2}
The operator $\wH_1$ is a bijective map when dealing with odd dimension, and its inverse
mapping is given by
\begin{equation}
\wH_1^{-1}(f)=\frac{\wH_m(1)}{m!}\sum_{j=0}^{(m-1)/2}(-1)^jc_{m,j}(m)\wH_1^{m-2j-1}(f).
\end{equation}
It is interesting to verify, by a straightforward calculation, that when the structural sets are identical, this operator reduces to the operator $\Psi^{-1}$ studied in \cite{morenodaniel}. This new reformulation of $\Psi^{-1}$, makes it possible to work directly with the full function $f$, without having to consider each of its $k$-vector parts separately.
\end{remark}
\begin{example}
Consider the following structural sets in $\R^3$:
\begin{equation*}
\varphi=\left\{\frac{\sqrt{2}}{2}e_1+\frac{\sqrt{2}}{2}e_2,\frac{\sqrt{2}}{2}e_1-\frac{\sqrt{2}}{2}e_2,-e_3\right\}
\end{equation*}
and
\begin{equation*}
\psi=\left\{e_3,e_2,e_1\right\}.
\end{equation*}
We want to find $f$ such that $\wH_1(f)=1$. In this case, we have
$$f=\frac{1}{6}[\wH_1^2(1)-7].$$
A straightforward computation shows that

\begin{align*}
\wH_1(1)&=\left(\frac{\sqrt{2}}{2}e_1+\frac{\sqrt{2}}{2}e_2\right)e_3+\left(\frac{\sqrt{2}}{2}e_1-\frac{\sqrt{2}}{2}e_2\right)e_2-e_3e_1\\
&=\frac{\sqrt{2}}{2}+\frac{\sqrt{2}}{2}e_1e_2+\left(\frac{\sqrt{2}}{2}+1\right)e_1e_3+\frac{\sqrt{2}}{2}e_2e_3,\\
\wH_1^2(1)&=\left(\frac{\sqrt{2}}{2}e_1+\frac{\sqrt{2}}{2}e_2\right)\left(\frac{\sqrt{2}}{2}+\frac{\sqrt{2}}{2}e_1e_2+\left(\frac{\sqrt{2}}{2}+1\right)e_1e_3+\frac{\sqrt{2}}{2}e_2e_3\right)e_3\\
&\quad+\left(\frac{\sqrt{2}}{2}e_1-\frac{\sqrt{2}}{2}e_2\right)\left(\frac{\sqrt{2}}{2}+\frac{\sqrt{2}}{2}e_1e_2+\left(\frac{\sqrt{2}}{2}+1\right)e_1e_3+\frac{\sqrt{2}}{2}e_2e_3\right)e_2\\
&\quad-e_3\left(\frac{\sqrt{2}}{2}+\frac{\sqrt{2}}{2}e_1e_2+\left(\frac{\sqrt{2}}{2}+1\right)e_1e_3+\frac{\sqrt{2}}{2}e_2e_3\right)e_1\\
&=\frac{1}{2}e_1e_3-\frac{1}{2}e_2e_3+\frac{1}{2}+\frac{\sqrt{2}}{2}-\frac{1}{2}e_1e_2+\frac{1}{2}e_2e_3+\frac{1}{2}e_1e_3+\left(\frac{1}{2}+\frac{\sqrt{2}}{2}\right)e_1e_2+\frac{1}{2}\\
&\quad+\frac{1}{2}e_1e_2+\frac{1}{2}+\left(\frac{1}{2}+\frac{\sqrt{2}}{2}\right)e_2e_3+\frac{1}{2}e_1e_3+\frac{1}{2}-\frac{1}{2}e_1e_2+\left(\frac{1}{2}+\frac{\sqrt{2}}{2}\right)e_1e_3-\frac{1}{2}e_2e_3\\
&\quad+\frac{\sqrt{2}}{2}e_1e_3+\frac{\sqrt{2}}{2}e_2e_3+\frac{\sqrt{2}}{2}+1+\frac{\sqrt{2}}{2}e_1e_2\\
&=3+\sqrt{2}+\sqrt{2}e_1e_2+(2+\sqrt{2})e_1e_3+\sqrt{2}e_2e_3,\\
\end{align*}
and, therefore,
$$f=-\frac{2}{3}+\frac{\sqrt{2}}{6}+\frac{\sqrt{2}}{6}e_1e_2+\left(\frac{1}{3}+\frac{\sqrt{2}}{6}\right)e_1e_3+\frac{\sqrt{2}}{6}e_2e_3.
$$
Now let us verify this:
\begin{align*}
\varphi_1f\psi_1&=\frac{\sqrt{2}}{6}+\frac{1}{3}+\frac{\sqrt{2}}{6}e_1e_2+\left(\frac{1}{3}-\frac{\sqrt{2}}{3}\right)e_1e_3-\frac{\sqrt{2}}{3}e_2e_3,\\
\varphi_2f\psi_2&=\frac{1}{3}-\frac{\sqrt{2}}{3}-\frac{\sqrt{2}}{3}e_1e_2+\left(\frac{1}{3}+\frac{\sqrt{2}}{6}\right)e_1e_3+\frac{\sqrt{2}}{6}e_2e_3,\\
\varphi_3f\psi_3&=\frac{1}{3}+\frac{\sqrt{2}}{6}+\frac{\sqrt{2}}{6}e_1e_2+\left(-\frac{2}{3}+\frac{\sqrt{2}}{6}\right)e_1e_3+\frac{\sqrt{2}}{6}e_2e_3.
\end{align*}
Thus,
$$\varphi_1f\psi_1+\varphi_2f\psi_2+\varphi_3f\psi_3=1.$$
Note that, in contrast to the case where the structural sets $\varphi$ and $\psi$ coincide, $f$ is not necessarily scalar-valued in the present setting.
\end{example}
When $m$ is odd, Proposition \ref{Prop2} (ii), together with the bijectivity of the operator $\wH_1$,  allows us to conclude that $f$ is $(\varphi,\psi)$-inframonogenic if and only if $\wH_1(f)$ is $(\varphi,\psi)$-inframonogenic. When $m$ is even, this can be established by means of the following proposition:
\begin{proposition}\label{Prop1dasdas}
Let $\varphi$ and $\psi$ be arbitrary structural sets, and let $m$ be even. 
If
\[
m\equiv 0 \pmod{4}
\quad\text{and}\quad
\det(\mathcal{M}_{\varphi,\psi})=1,
\]
or
\[
m\equiv 2 \pmod{4}
\quad\text{and}\quad
\det(\mathcal{M}_{\varphi,\psi})=-1,
\]
then an $\R_{0,m}^{-}$-valued function $f$ is $(\varphi,\psi)$-inframonogenic if and only if $\wH_1(f)$ is $(\varphi,\psi)$-inframonogenic. Conversely, if
\[
m\equiv 2 \pmod{4}
\quad\text{and}\quad
\det(\mathcal{M}_{\varphi,\psi})=1,
\]
or
\[
m\equiv 0 \pmod{4}
\quad\text{and}\quad
\det(\mathcal{M}_{\varphi,\psi})=-1,
\]
then an $\R_{0,m}^{+}$-valued function $f$ is $(\varphi,\psi)$-inframonogenic if and only if $\wH_1(f)$ is $(\varphi,\psi)$-inframonogenic.
\end{proposition}

Let us now consider the equation
\begin{equation}\label{eqseccfinal}
df+\wH_1^2(f)=g,\quad d\in\R.
\end{equation}
In this case, if $m=2l$ we have
\begin{equation}
\wH_{2l}(f)=\frac{1}{(2l)!}\sum_{j=0}^l(-1)^jc_{2l,j}(m)\left[\sum_{r=0}^{l-j-1}(-d)^r\wH_1^{2l-2j-2r-2}(g)+(-d)^{l-j}f\right],
\end{equation} 
and, therefore,
\begin{align}
&\left[1-\frac{\wH_m(1)}{(2j)!}\sum_{j=0}^l(-1)c_{2l,j}(m)(-d)^{l-j}\right]f_+-\left[1+\frac{\wH_m(1)}{(2j)!}\sum_{j=0}^l(-1)c_{2l,j}(m)(-d)^{l-j}\right]f_-\nonumber\\&=\frac{\wH_m(1)}{(2l)!}\sum_{j=0}^l(-1)^jc_{2l,j}(m)\sum_{r=0}^{l-j-1}(-d)^r\wH_1^{2l-2j-2r-2}(g).\label{forpr3}
\end{align}
If $m=2l+1$, then we obtain
\begin{equation}
\wH_{2l+1}(f)=\frac{1}{(2l+1)!}\sum_{j=0}^l(-1)^jc_{2l+1,j}(m)\left[\sum_{r=0}^{l-j-1}(-d)^r\wH_1^{2l-2j-2r-1}(g)+(-d)^{l-j}\wH_1(f)\right],
\end{equation} 
and, thus,
\begin{align}
&f+\frac{\wH_m(1)}{(2l+1)!}\sum_{j=0}^l(-1)^jc_{2l+1,j}(m)(-d)^{l-j}\wH_1(f)\nonumber\\&=-\frac{\wH_{m}(1)}{(2l+1)!}\sum_{j=0}^l(-1)^jc_{2l+1,j}(m)\sum_{r=0}^{l-j-1}(-d)^r\wH_1^{2l-2j-2r-1}(g).\label{forpr}
\end{align}
Applying the operator $\wH_1$, we obtain
\begin{align}
&\wH_1(f)+\frac{\wH_m(1)}{(2l+1)!}\sum_{j=0}^l(-1)^jc_{2l+1,j}(m)(-d)^{l-j}(g-df)\nonumber\\
&=-\frac{\wH_{m}(1)}{(2l+1)!}\sum_{j=0}^l(-1)^jc_{2l+1,j}(m)\sum_{r=0}^{l-j-1}(-d)^r\wH_1^{2l-2j-2r}(g).\label{forpr2}
\end{align}
The formulas \eqref{forpr3}, \eqref{forpr}, and \eqref{forpr2} allow us to explicitly determine $f$ in terms of $g$ and will be useful in the next section.

\section{Fischer decompositions}
An inner product on the space $\cP(k)$ of $\R_{0,m}$-valued $k$-homogeneous
polynomials is defined by
\begin{equation}\label{product}
\langle P_k,Q_k\rangle
=
\big[\overline{P_k(\pux)}\,Q_k\big]_0,
\qquad
P_k,Q_k\in\cP(k).
\end{equation}
This is the standard Fischer inner product in the Clifford-algebraic
setting.

Let $\varphi=\{\varphi_1,\ldots,\varphi_m\}$ and
$\psi=\{\psi_1,\ldots,\psi_m\}$ be arbitrary structural sets. 
The following adjointness relations with respect to the Fischer inner
product will be used repeatedly. If
$P_{k-1}\in\cP(k-1)$ and $Q_k\in\cP(k)$, then
\begin{align*}
\langle \ux_\varphi P_{k-1},Q_k\rangle
&=
\big[\overline{\hD P_{k-1}(\pux)}\,Q_k\big]_0
\nonumber\\
&=
-\big[\overline{P_{k-1}(\pux)}\,\hD Q_k\big]_0
=
-\langle P{k-1},\hD Q_k\rangle,
\\[1ex]
\langle P_{k-1}\ux_\psi,Q_k\rangle
&=
\big[\overline{P_{k-1}(\pux)\pD}\,Q_k\big]_0
\nonumber\\
&=
-\big[\overline{P_{k-1}(\pux)}\,Q_k\pD\big]_0
=
-\langle P_{k-1},Q_k\pD\rangle.
\end{align*}
Thus, multiplication from the left by $\ux_\varphi$ and from the right by
$\ux_\psi$ are, respectively, adjoint to the differential operators
$-\hD$ and $-\pD$ with respect to the Fischer inner product.

Consequently, for
$P_{k-2}\in\cP(k-2)$ and $Q_k\in\cP(k)$, we obtain
\begin{align}
\label{f1}
\langle \ux_\varphi P_{k-2}\ux_\psi,Q_k\rangle
&=
-\langle P_{k-2}\ux_\psi,\hD Q_k\rangle
=
\langle P_{k-2},\hD Q_k\pD\rangle,
\\[1ex]
\label{f2}
\langle \ux_\psi \ux_\varphi P_{k-2},Q_k\rangle
&=
-\langle \ux_\varphi P_{k-2},\pD Q_k\rangle
=
\langle P_{k-2},\hD\pD Q_k\rangle.
\end{align}
These identities will play a fundamental role in the Fischer
decompositions considered below. Let $\I_{\varphi,\psi}(k)\subset\cP(k)$ denotes the set of all $\R_{0,m}$-valued $(\varphi,\psi)$-inframonogenic $k$-homogeneous polynomials in $\R^m$. The following theorem was established in \cite{FischerDAS}:
\begin{theorem}{\cite[Theorem 3.1, p. 4]{FischerDAS}}\label{FD1}
Let be $k\geq 2$, then the following decomposition holds:
\[
\cP(k)=\I_{\varphi,\psi}(k)\oplus \ux_\varphi\cP(k-2)\ux_{\psi}.
\]
Moreover, the subspaces $\I_{\varphi,\psi}(k)$ and $\ux_\varphi\cP(k-2)\ux_{\psi}$ are orthogonal w.r.t the inner product \eqref{product}.
\end{theorem}
Thanks to the preceding theorem, we can derive the complete Fischer decomposition in terms of $(\varphi,\psi)$-inframonogenic polynomials mentioned in the Introduction:
\begin{equation}\label{fischerinfrageneral2}
\cP(\R^m,\R_{0,m})=\bigoplus_{k=0}^{\infty}\bigoplus_{p=0}^{\infty}\ux_\varphi^p\I_{\varphi,\psi}(k)\ux_\psi^p,
\end{equation}
Since $$\hD\hD(.)\pD\pD=\Delta_m^2,$$ the class of $(\varphi,\psi)$-inframonogenic polynomials forms a proper subspace of the space of biharmonic polynomials. The following theorem is a classical and well-known result in real and Clifford analysis. Although its details of the proof can be found in standard references on these subjects \cite{DSS}, we include it here for the reader's convenience and to make the discussion self-contained.

\begin{theorem}[Biharmonic Fischer decomposition]\label{FD1biar}
Let be $k\geq 4$, then the following decomposition holds:
\[
\cP(k)=\B(k)\oplus |\ux|^4\cP(k-4),
\]
where $\B(k)$ denotes the subspace of all $\R_{0,m}$-valued biharmonic $k$-homogeneous polynomials in $\R^m$. Moreover, the subspaces $\B(k)$ and $|\ux|^4\cP(k-4)$ are orthogonal w.r.t the inner product \eqref{product}.
\end{theorem}

\pf Since $\cP(k)=|\ux|^4P(k-4)\oplus (|\ux|^4P(k-4))^{\perp}$, we can restrict to prove that 
\begin{equation}\label{relbiha}
\B(k)=(|\ux|^4P(k-4))^{\perp}
\end{equation}
Indeed, let be $P_k\in (|\ux|^4P(k-4))^{\perp}$. Then for any $Q_{k-4}\in\cP(k-4)$ one has 
\[\langle |\ux|^4 Q_{k-4}, P_k\rangle
=\langle \ux^4 Q_{k-4}, P_k\rangle=\langle Q_{k-4},\pux^4 P_k\rangle=\langle Q_{k-4},\Delta_m^2 P_k\rangle=0,
\]
where use has been made of \eqref{f2} with $\varphi=\psi=\{e_1,...,e_m\}$. In particular, for $Q_{k-4}=\Delta_m^2 P_k$, we arrive at $\Delta_m^2 P_k=0$ and so $P_k\in\B(k)$. Consequently, 
\[
(|\ux|^4\cP(k-4))^{\perp}\subset\B(k).
\] 

Now let $P_k\in\B(k)$. Then, for any $Q_{k-4}\in\cP(k-4)$ we also have
\[
\langle |\ux|^4 Q_{k-4}, P_k\rangle=\langle Q_{k-4},\Delta_m P_k\rangle=0,
\]
and so $P_k\in (|\ux|^4\cP(k-4))^{\perp}$.

From this we deduce that $\B(k)\subset (|\ux|^4\cP(k-4))^{\perp}$ and, finally, relation \eqref{relbiha} follows. 
\qed
The following result is an immediate consequence of Theorem \ref{FD1}.
\begin{theorem}[Complete biharmonic Fischer decomposition]\label{FD2}
Let be $k\geq 4$, then the following decomposition holds:
\begin{equation}\label{biharmonicFischer}
\cP(\R_m,\R_{0,m})=\bigoplus_{k=0}^{\infty}\bigoplus_{p=0}^{\infty}|\ux|^{4p}\B(k).
\end{equation}
\end{theorem}

To summarize the preceding results, let us present a list of Fischer-type decompositions of the space of $\R_{0,m}$-valued $k$-homogeneous polynomials in $\R^m$:
\begin{align}
\mathcal{P}(k)&=\bigoplus_{s=0}^{ k}\ux^s\mathcal{M}(k-s),\label{monogenic}\\
\mathcal{P}(k)&=\bigoplus_{s=0}^{\lfloor k/2\rfloor}\ux^{2s}\mathcal{H}(k-2s),\label{harmonic}\\
\mathcal{P}(k)&=\bigoplus_{s=0}^{\lfloor k/2\rfloor}\ux^s\mathcal{I}(k-2s)\ux^s,\label{inframonogenic}\\
\mathcal{P}(k)&=\bigoplus_{s=0}^{\lfloor k/2\rfloor}\ux_\varphi^{s}\mathcal{I}_{\varphi,\psi}(k-2s)\ux_\psi^s,\label{inframonogenicgen}\\
\mathcal{P}(k)&=\bigoplus_{s=0}^{\lfloor k/4\rfloor}\ux^{4s}\mathcal{B}(k-4s).\label{biharmonic}
\end{align}
The first thing we would like to point out is that decomposition \eqref{monogenic} can be viewed as a refinement of decomposition \eqref{harmonic}. Indeed, for every $k$-homogeneous polynomial $P_k$, there exist uniquely determined left-monogenic $j$-homogeneous polynomials $M_j$, $j=0,...,k$, such that
\begin{equation}\label{auxiliar}
P_k=M_k+\ux M_{k-1}+\ux^2(M_{k-2}+\ux M_{k-3})+...+\ux^{2s}(M_{k-2s}+\ux M_{k-2s-1})+...+\ux^kM_0,
\end{equation}
which follows from \eqref{monogenic}.

Moreover, observe that all $k$-homogeneous polynomials of the form $M_{k}+\ux M_{k-1}$ are harmonic: 
\begin{align*}
\pux[M_{k}+\ux M_{k-1}]&=-(2k-2+m)M_{k-1},\\
\Delta_m[M_{k}+\ux M_{k-1}]&=(2k-2+m)\pux M_{k-1}=0.
\end{align*}
Consequently, decomposition \eqref{auxiliar} is nothing but the Fischer decomposition of $k$-homogeneous polynomials in terms of harmonic polynomials.

In particular, by the uniqueness of decomposition \eqref{harmonic}, every harmonic $k$-homogeneous polynomial $H_k$ can be decomposed as
\begin{equation}\label{monogenicrelhar}
H_k = M_k + \ux M_{k-1}.
\end{equation}

Similarly, decomposition \eqref{harmonic} can also be viewed as a refinement of \eqref{biharmonic}. By virtue of decomposition \eqref{harmonic}, for every $k$-homogeneous polynomial $P_k$, there exist uniquely determined harmonic $j$-homogeneous polynomials $H_j$, $j=k-2\lfloor l/2\rfloor,...,k$, such that 
\begin{equation}\label{42ref}
P_k=H_k+\ux^2H_{k-2}+\ux^4(H_{k-4}+\ux^2H_{k-6})+...+\ux^{4s}(H_{k-4s}+\ux^2H_{k-4s-2})+...+\ux^{2\lfloor k/2\rfloor}H_{k-2\lfloor k/2\rfloor}.
\end{equation}
In this case, note that the $k$-homogeneous polynomials $H_{k}+\ux^{2}H_{k-2}$, are biharmonic:
\begin{align*}
\pux [H_{k}+\ux^2 H_{k-2}]&=\pux H_k-2\ux H_{k-2}+\ux^2(\pux H_{k-2}),\\
\pux^2 [H_k+\ux^2 H_{k-2}]&=2(2k-4+m)H_{k-2}+2\ux(\pux H_{k-2})-2\ux(\pux H_{k-2}),\\
\Delta_m^2 [H_k+\ux^2 H_{k-2}]&=2(2k-4+m)\Delta_m H_{k-2}=0.
\end{align*}
Thus, if $P_k$ is biharmonic, then the uniqueness of the decomposition \eqref{biharmonic} allows us to identify the corres\-ponding biharmonic components in \eqref{42ref}, thereby yielding relation
\begin{equation}\label{relbiharmonicharmonic}
B_k=H_k+\ux^2 H_{k-2},
\end{equation}
where $B_k$ is any biharmonic $k$-homogeneous polynomial. Using relation \eqref{monogenicrelhar} and analogous arguments, we can also obtain
\begin{equation}
B_k=M_k+\ux M_{k-1}+\ux^2M_{k-2}+\ux^3M_{k-3},
\end{equation}
where $M_j$, $j\in\{k,k-1,k-2,k-3\}$, are left-monogenic $j$-homogeneous polynomials.

Special attention should be given to decompositions \eqref{inframonogenic} and \eqref{inframonogenicgen}. Since the spaces of inframonogenic functions and, more generally, $(\varphi,\psi)$-inframonogenic functions, are subspaces of the space of biharmonic functions, it is natural to ask whether these decompositions can be regarded as refinements of the Fischer decomposition in terms of biharmonic polynomials

By \eqref{inframonogenicgen}, and using the fact that $\ux_\varphi^2=\ux_\psi^2=\ux^2$, we obtain
\begin{equation}\label{45ref}
P_k=I_k+\ux_{\varphi}I_{k-2}\ux_{\psi}+\ux^4(I_{k-4}+\ux_{\varphi}I_{k-6}\ux_\psi)+...+\ux^{4s}(I_{k-4s}+\ux_{\varphi}I_{k-4s-2}\ux_\psi)+...+\ux_\varphi^{\lfloor k/2\rfloor}I_{k-2\lfloor k/2\rfloor}\ux_\psi^{\lfloor k/2\rfloor},
\end{equation}
where $I_j$ denotes $(\varphi,\psi)$-inframonogenic $j$-homogeneous  polynomials.

If decomposition \eqref{45ref}  were a refinement of Fischer decomposition by biharmonic polynomials, one would naturally expect the polynomials $I_k+\ux_{\varphi}I_{k-2}\ux_\psi$ to be biharmonic. More precisely, we are interested in determining whether every $B_k\in\mathcal{B}(k)$
admits a unique decomposition of the form
\begin{equation}\label{birelproof}
B_k=I_k+\ux_{\varphi}I_{k-2}\ux_{\psi},
\end{equation}
where $I_k$ and $I_{k-2}$ are $(\varphi,\psi)$-inframonogenic homogeneous polynomials of degrees $k$ and $k-2$, respectively. As we shall see, the answer is negative.

Indeed, we have
{\footnotesize{
\begin{align*}
\hD[\underline{x}_\varphi I_{k-2}\underline{x}_\psi]&=(\hD[\underline{x}_\varphi I_{k-2}])\underline{x}_\psi+{\wH_1}(\underline{x}_\varphi I_{k-2}),\\
\hD\hD[\underline{x}_\varphi I_{k-2}\underline{x}_\psi]&=(\hD\hD[\underline{x}_\varphi I_{k-2}])\underline{x}_\psi-2[\underline{x}_\varphi I_{k-2}]\pD,\\
\hD\hD\hD[\underline{x}_\varphi I_{k-2}\underline{x}_\psi]&=(\hD\hD\hD[\underline{x}_\varphi I_{k-2}])\underline{x}_\psi+{\wH_1}(\hD\hD[\underline{x}_\varphi I_{k-2}])-2\hD[\underline{x}_\varphi I_{k-2}]\pD,\\
\Delta_m^2[\underline{x}_\varphi I_{k-2}\underline{x}_\psi]
&=\underline{x}_\varphi(\Delta_m^2I_{k-2})\underline{x}_\psi+8\hD[I_{k-2}]\pD+4{\wH_1}(\Delta_m I_{k-2})+4\underline{x}_\varphi(\Delta_m[I_{k-2}]\pD)+4(\hD\Delta_m[I_{k-2}])\underline{x}_\psi\\&=4{\wH_1}(\Delta_m I_{k-2}).
\end{align*}}}
Therefore, if \eqref{birelproof} holds, then the
$(\varphi,\psi)$-inframonogenic polynomial $I_{k-2}$ must be either harmonic or such that $\varpi I_{k-2}$ is
$\R_{0,m}^{(m/2)}$-valued when $m$ is even, where $\varpi$ is an
invertible element of the Clifford algebra satisfying
\[
\varpi\varphi_j=\psi_j\varpi,
\qquad j=1,\ldots,m.
\]

Let $P_4$ be a polynomial in $\cP(4)$. By \eqref{inframonogenicgen}, we have
$$P_4=I_4+\ux_\varphi I_2\ux_\psi+\ux_\varphi^2 I_0\ux_\psi^2,$$
where $I_0$, $I_2$ and $I_4$ are $(\varphi,\psi)$-inframonogenic homogeneous polynomials of degrees $0$, $2$ and $4$, respectively.  Using Propositions \ref{Prop1} and \ref{Prop2}, a direct computation yields
\begin{align*}
\hD P_4&=\hD I_4-(6+m)I_2\ux_\psi-\ux_\varphi\hD[I_2\ux_\psi]-2\ux_\varphi I_0\ux_\psi^2+\ux_\varphi^2\hD[I_0\ux_\psi^2]\\
&=\hD I_4-(6+m)I_2\ux_\psi-\ux_\varphi(\hD I_2)\ux_\psi-\ux_\varphi\wH_1(I_2)-2\ux_\varphi I_0\ux_\psi^2-2\ux_\varphi^3I_0,\\
\hD P_4\pD&=(6+m)(4+m)I_2+(6+m)(I_2\pD)\ux_\psi+(4+m)\ux_\varphi(\hD I_2)+\wH_1(\hD I_2)\ux_\psi\\
&\quad+2\ux_\varphi (\hD I_2)+\ux_\varphi\wH_1(I_2\pD)-\wH_1^2(I_2)+4\ux_\varphi I_0\ux_\psi-2\wH_1(I_0)\ux_\psi^2\\
&\quad-2\ux_\varphi^2\wH_1(I_0)+4\ux_\varphi I_0\ux_\psi\\
&=(6+m)(4+m)I_2-\wH_1^2(I_2)+(6+m)(I_2\pD)\ux_\psi+(6+m)\ux_\varphi(\hD I_2)\\
&\quad+\wH_1(\hD I_2)\ux_\psi+\ux_\varphi\wH_1(I_2\pD)-4\ux_\varphi^2\wH_1(I_0)+8\ux_\varphi I_0\ux_\psi,\\
\hD\!\!^2 P_4\pD&=(6+m)(4+m)\hD I_2-4\hD I_2-4\wH_1(I_2\pD)-\wH_1^2(\hD I_2)+(6+m)\wH_1(I_2\pD)\\
&\quad-(6+m)(2+m)\hD I_2-(6+m)\ux_\varphi(\hD\!\!^2I_2)-\wH_1(\hD\!\!^2 I_2)\ux_\psi+\wH_1^2(\hD I_2)\\&\quad-(2+m)\wH_1(I_2\pD)+2\ux_\varphi(I_2\pD\!\!^2)+8\ux_\varphi\wH_1(I_0)-8(2+m)I_0\ux_\psi-8\ux_\varphi\wH_1(I_0)\\
&=2(m+4)\hD I_2-(4+m)\ux_\varphi(\hD\!\!^2 I_2)-\wH_1(\hD\!\!^2 I_2)\ux_\psi-8(2+m)I_0\ux_\psi,\\
\hD P_4\pD\!\!^2&=2(m+4)I_2\pD-(4+m)(I_2\pD\!\!^2)\ux_\psi-\ux_\varphi\wH_1(I_2\pD\!\!^2)-8(2+m)\ux_\varphi I_0,\\
\hD\!\!^2 P_4\pD\!\!^2&=-4\wH_1(\hD\!\!^2 I_2)+8m(2+m)I_0,\\
\hD\!\!^3P_4\pD&=(2+m)(m+4)\hD\!\!^2I_2-\wH_1^2(\hD\!\!^2I_2)-8(2+m)\wH_1(I_0).
\end{align*}
Therefore,
\begin{equation}\label{45}
-m\Delta_m\hD P_4\pD+\wH_1(\Delta_m^2 P_4)=-m(2+m)(4+m)\Delta_m I_2+(m+4)\wH_1^2(\Delta_m I_2).
\end{equation}
Note that equation \eqref{45} has a form similar to that of \eqref{eqseccfinal}, which was analyzed in the previous section. Formulas \eqref{forpr3}, \eqref{forpr}, and \eqref{forpr2} allow us to express $\Delta_m I_2$ solely in terms of $\Delta_m \hD P_4\pD$ and $\Delta_m^2 P_4$. From the relations established above, we obtain
\[
I_0=\frac{1}{8m(2+m)}
\left[\Delta_m^2P_4-4\wH_1(\Delta_m I_2)\right],
\]
where, in this case, $\Delta_m I_2$ is expressed as a sum involving
$\Delta_m\hD P_4\pD$ and $\Delta_m^2P_4$. Therefore, in general, if $P_4$ is a biharmonic polynomial, then $I_0$ does not necessarily vanish, since $\Delta_m\hD P_4\pD$ need not vanish. This explains why, in general, a relation such as \eqref{birelproof} does not hold. For example, when $\varphi=\psi=\{e_1,...,e_m\}$ we obtain
\begin{align*}
I_0&=\frac{1}{8m(m+2)}\Delta_m^2 P_4\\&\quad+\sum_{k=0}^m\frac{1}{4m(m+2)(m+4)(m+2km-2k^2)}(m(-1)^{k+1}(2k-m)[\Delta_m\pux P_4\pux]_k+(2k-m)^2[\Delta_m^2 P_4]_k).
\end{align*}
Furthermore, after some additional computations, we also obtain that

\begin{align*}
\pux I_2&=\frac{2}{m+4}\pux^2 P_4\pux+\frac{1}{2}\ux \sum_{k=0}^m\frac{1}{2(m+4)(m+2km-2k^2)}(m[\pux^3 P_4\pux]_k+(-1)^k(2k-m)[\pux^2 P_4\pux^2]_k)\\
&\quad+\sum_{k=0}^m\frac{(-1)^k(2k-m)}{4(m+4)^2(m+2km-2k^2)}(m[\pux^3 P_4\pux]_k+(-1)^k(2k-m)[\pux^2 P_4\pux^2]_k)\ux\\&\quad+\frac{4(m+2)}{m+4}\frac{1}{8m(m+2)}(\pux^2 P_4\pux^2)\ux\\&\quad+\sum_{k=0}^m\frac{1}{4m(m+2)(m+4)(m+2km-2k^2)}(m(-1)^k(2k-m)[\pux^3 P_4\pux]_k+(2k-m)^2[\pux^2 P_4\pux^2]_k)\ux,\\
I_2\pux &=\frac{2}{m+4}\pux P_4\pux^2+\frac{1}{2}\sum_{k=0}^m\frac{1}{2(m+4)(m+2km-2k^2)}(m[\pux^3 P_4\pux]_k+(-1)^k(2k-m)[\pux^2 P_4\pux^2]_k)\ux\\
&\quad+\sum_{k=0}^m\frac{(-1)^k(2k-m)}{4(m+4)^2(m+2km-2k^2)}\ux(m[\pux^3 P_4\pux]_k+(-1)^k(2k-m)[\pux^2 P_4\pux^2]_k)\\&\quad+\frac{4(m+2)}{m+4}\frac{1}{8m(m+2)}\ux(\pux^2 P_4\pux^2)\\&\quad+\sum_{k=0}^m\frac{1}{4m(m+2)(m+4)(m+2km-2k^2)}\ux(m(-1)^k(2k-m)[\pux^3 P_4\pux]_k+(2k-m)^2[\pux^2 P_4\pux^2]_k),\\
I_2&=\sum_{k=0}^m\frac{1}{2(5m+2km-2k^2+12)}[\pux P_4\pux-(6+m)(I_2\pux)\ux+(6+m)\ux(\pux I_2)+\wH_1[\pux I_2]\ux\\&\quad+\ux\wH_1[I_2\pux]+8\ux I_0\ux-4x^2\wH_1[I_0]]_k,\\
I_4&=P_4-\ux I_2\ux-\ux^2 I_0\ux^2. 
\end{align*}
Therefore, this shows that, in general, a biharmonic $4$-homogeneous polynomial $B_4$ admits a Fischer decomposition \eqref{inframonogenicgen} involving all the terms:
\begin{equation}
B_4=I_4+\ux_\varphi I_2\ux_\psi+\ux_\varphi^2 I_0\ux_\psi^2.
\end{equation}
More generally, every biharmonic $k$-homogeneous polynomial $B_k$ admits a complete Fischer decomposition of the form \eqref{inframonogenicgen}, in which all the summands contribute to the decomposition.

This further highlights the essential distinction between the Fischer decompositions \eqref{inframonogenic} and \eqref{inframonogenicgen} and the traditional Fischer decomposition by biharmonic polynomials in \eqref{biharmonic}, showing that the former cannot, in general, be regarded as a conventional refinement of the latter.	

This observation reveals a fundamental distinction between the Fischer decomposition in terms of inframonogenic polynomials, and more generally, $(\varphi,\psi)$-inframonogenic polynomials, and the more classical decomposition based on biharmonic polynomials. In the latter setting, the biharmonic component is determined solely by the underlying structure of Laplacian, whereas in the $(\varphi,\psi)$-inframonogenic setting, the decomposition is influenced by the choice of the structural sets $\varphi$ and $\psi$. This additional freedom makes it possible to obtain different classes of polynomial decompositions while preserving the under\-lying Clifford-algebraic structure. Consequently, the appropriate choice of structural sets gives rise to a broad family of Fischer-type decompositions, providing new ways of representing polynomial spaces and exten\-ding the traditional theory beyond the biharmonic framework. 

A similar phenomenon occurs for the class of $(\varphi,\psi)$-harmonic polynomials, namely, the solutions of the equation
\begin{equation}
\hD\pD f=0,
\end{equation}
for which one also obtains a complete Fischer decomposition of the form:
\begin{equation}
\mathcal{P}(k)=\bigoplus_{s=0}^{\lfloor k/2\rfloor}(\ux_\psi\ux_\varphi)^{s}\mathcal{H}_{\varphi,\psi}(k-2s),
\end{equation}
where $\mathcal{H}_{\varphi,\psi}(k)$ denotes the subspace of all $\R_{0,m}$-valued $(\varphi,\psi)$-harmonic $k$-homogeneous polynomials in $\R^m$ (see \cite[Theorem 5.1, p. 11]{FischerDAS}). 

It is important to observe that both the class of $(\varphi,\psi)$-harmonic polynomials and the class of $(\varphi,\psi)$-inframonogenic polynomials constitute nontrivial subclasses of the space of biharmonic polynomials \cite{S1}. These inclusions provide a natural bridge between the corresponding Fischer decompositions and the theory of biharmonic polynomials, while revealing additional structures induced by the choice of the structural sets $\varphi$ and $\psi$. In particular, although every polynomial belonging to either of these classes is biharmonic, the converse does not hold in general. Thus, the introduction of the structural sets $\varphi$ and $\psi$ provides a richer framework in which one can distinguish and study special subclasses within the broader family of biharmonic polynomials. In this case, we likewise cannot guarantee that every polynomial $B_k\in\mathcal{B}(k)$ admits a decomposition of the form
\begin{equation}
B_k=G_k+\ux_\psi\ux_\varphi G_{k-2},
\end{equation}
where $G_k$ and $G_{k-2}$ are $(\varphi,\psi)$-harmonic homogeneous polynomials of degrees $k$ and $k-2$, respectively. Indeed, calculations and arguments analogous to those presented above lead to the same conclusion.
\section*{Concluding remarks}
The study of $(\varphi,\psi)$-inframonogenic and $(\varphi,\psi)$-harmonic functions reveals a rich class of Clifford algebra-valued functions that extends several classical settings in Clifford analysis. Both classes are defined through second-order differential equations involving two structural sets, providing a flexible framework in which the choice of $\varphi$ and $\psi$ influences the underlying differential structure and the resulting decomposition properties.

A particularly relevant feature is that both $(\varphi,\psi)$-inframonogenic and $(\varphi,\psi)$-harmonic polynomials form nontrivial subclasses of the space of biharmonic polynomials. Thus, although they are closely related to the traditional biharmonic theory, they retain additional algebraic information determined by the structural sets. This distinction becomes especially apparent in the corresponding Fischer decompositions, where the components are described in terms of $(\varphi,\psi)$-harmonic or $(\varphi,\psi)$-inframonogenic polynomials rather than solely in terms of the harmonic components.

It is worth emphasizing that the classical Fischer decomposition by harmonic polynomials is recovered, in the $(\varphi,\psi)$-harmonic setting, when the two structural sets coincide, that is, when $\varphi=\psi$. In contrast, the $(\varphi,\psi)$-inframonogenic decomposition should not, in general, be regarded as a mere reformulation of the harmonic Fischer decomposition, even when the structural sets coincide. This distinction reveals the genuinely different nature of the inframonogenic decomposition.

The freedom in choosing the structural sets $\varphi$ and $\psi$ therefore gives rise to a broad family of Fischer-type decompositions. In particular, the $(\varphi,\psi)$-harmonic framework contains the classical harmonic decomposition as a special case, while more general choices of structural sets produce new representations that reflect the underlying Clifford-algebraic structure. This flexibility makes both classes a natural setting for further investigations of generalized Fischer decompositions, boundary value problems, and related Dirac-type systems.

It is also worth emphasizing that, in this work, we obtained an explicit formula for the inverse of the operator $\wH_1$, which is of particular relevance both to the present study and to future developments. Interestingly, the formula stated in Remark \ref{Remark3.2} reduces, in the standard scenario, to the one obtained in \cite{morenodaniel}. Given the variety of structural sets involved, and the fact that the operator $\wH_1$ does not necessarily preserve $k$-vectors, an explicit formula for its inverse, whenever it exists, is of considerable interest in this setting.
\section*{Declarations}
\subsection*{Author contribution statement}
All authors contributed equally to this work
\subsection*{Funding}
Postgraduate Study Fellowship of the Secretar\'ia de Ciencia, Humanidades, Tecnolog\'ia e Innovaci\'on  (grant number 1043969).
\subsection*{Conflict of interest}
The authors declare that they have no conflict of interest regarding the publication of this paper.
\subsection*{Data Availability Statement}
Data sharing is not applicable to this article as no datasets were generated or analyzed during the current study.

\end{document}